# Modified discrete random walk with absorption

**Theo van Uem**

Amsterdam School of Technology, Weesperzijde 190, 1097 DZ Amsterdam, The Netherlands.

Email: t.j.van.uem@hva.nl

Abstract

We obtain expected number of arrivals, probability of arrival, absorption probabilities and expected time before absorption for a modified discrete random walk on the (sub)set of integers. In a [pqrs] random walk the particle can move one step forward or backward, stay for a moment in the same state or it can be absorbed immediately in the current state. M[pqrs] is a modified version, where probabilities on both sides of a multiple function barrier M are of different [pqrs] type.

## 1. Introduction

Random walk can be used in various disciplines: in physics as a simplified model of Brownian motion, in ecology to describe individual animal movements and population dynamics, in statistics to analyze sequential test procedures, in economics to model share prices and their derivatives, in medicine and biology where absorbing barriers give a natural model for a wide variety of phenomena. Random walks have been studied for decades on regular structures such as lattices. We now give a brief historical review of the use of barriers in a one-dimensional discrete random walk. Weesakul (1961) discussed the classical problem of random walk restricted between a reflecting and an absorbing barrier. Using generating functions he obtains explicit expressions for the probability of absorption. Lehner (1963) studies one-dimensional random walk with a partially reflecting barrier using combinatorial methods. Gupta (1966) introduces the concept of a multiple function barrier (mfb): a state that can absorb, reflect, let through or hold for a moment. Dua, Khadilkar and Sen (1976) find the bivariate generating functions of the probabilities of a particle reaching a certain state under different conditions. Percus (1985) considers an asymmetric random walk, with one or two boundaries, on a one-dimensional lattice. At the boundaries, the walker is either absorbed or reflected back to the system. Using generating functions the probability distribution of being at position m after n steps is obtained, as well as the mean number of steps before absorption. El-Shehawey (2000) obtains absorption probabilities at the boundaries for a random walk between one or two partially absorbing boundaries as well as the conditional mean for the number of steps before stopping given the absorption at a specified barrier, using conditional probabilities.

In this paper we investigate one dimensional random walk where every state can be considered as a multiple function barrier (mfb). A mfb can absorb, reflect, let through or hold for a moment. The difference with earlier investigations is that the number of barriers is not limited: every state is a mfb. A [pqrs] random walk on the integers is defined by: p is the one-step forward probability, q one-step backward, r the probability to stay for a moment in the same position and s is the probability of immediate absorption in the current state (p+q+r+s=1, pqs>0).

A M[pqrs] random walk on the integers is a modified [pqrs] random walk, where probabilities right and left of the mfb M are different: $[p_1, q_1, r_1, s_1]$ and $[p_2, q_2, r_2, s_2]$.

In section 2 we solve sets of difference equations which are related to the expected number of arrivals and expected time before absorption.

In sections 3-5 we study [pqrs] random walk on $[0, N], [0, \infty)$ and $(-\infty, \infty)$.

We do the same for M[pqrs] random walk in sections 6-8.

## 2. Sets of difference equations.

We define a probability generating function:

$$X_j = X_j(z) = X_{i,j}(z) = \sum_{k=0}^{\infty} p_{ij}^{(k)} z^k \qquad (0 < z \leq 1)$$

The expected number of visits to state j is given by:

$$x_j = X_j(1) = X_{i,j}(1) = \sum_{k=0}^{\infty} p_{ij}^{(k)}$$

We define:

$m_i$ = expected time before absorption when starting in i.

$m_{ij}$ = expected time before absorption in j when starting in i.

The following two theorems are the base of the results in sections 3-8.

**Theorem 1**

The set of difference equations:

$$(1-rz)X_n = \delta(n, i_0) + pzX_{n-1} + qzX_{n+1} \quad (n \in \mathbb{Z}) \qquad (pq > 0,\ p+q+r < 1,\ 0 < z \leq 1) \qquad (1)$$

has solutions:

$$X_n(z) = \begin{cases} \zeta \xi_1^{n-i_0} + C_1 \xi_1^n + C_2 \xi_2^n & (n \leq i_0) \\ \zeta \xi_2^{n-i_0} + C_1 \xi_1^n + C_2 \xi_2^n & (n \geq i_0) \end{cases}$$

where:

$$\xi_1(z) = \frac{(1-rz) + \sqrt{(1-rz)^2 - 4pqz^2}}{2qz} > 1$$

$$0 < \xi_2(z) = \frac{(1-rz) - \sqrt{(1-rz)^2 - 4pqz^2}}{2qz} < 1$$

$$\zeta = \zeta_z = [(1-rz)^2 - 4pqz^2]^{-\frac{1}{2}}$$

Proof:

General solution of homogeneous part of (1) is:

$$X_n = C_1 \xi_1^n + C_2 \xi_2^n \quad (n \in \mathbb{Z}) \text{, where } qz\xi^2 - (1-rz)\xi + pz = 0$$

We find a particular solution of (1) by applying the method of variation of the constants in:

$$X_n = \begin{cases} C_1 \xi_1^n & (n \leq i_0) \\ C_2 \xi_2^n & (n \geq i_0) \end{cases}$$

Note:

We will use $\xi = \xi(z)$ in combination with X. $\xi = \xi(1)$ will be used in combination with x and m.

**Theorem 2**

The set of difference equations:

$$(1-r)m_i = pm_{i+1} + qm_{i-1} + 1 - s \quad (i \in \mathbb{Z}) \ (p+q+r+s = 1,\ pqs > 0)$$

has solutions

$$m_i = a\xi_1^{-i} + b\xi_2^{-i} + \frac{1-s}{s} \quad (i \in \mathbb{Z})$$

Proof:

By substitution.

## 3. [pqrs] random walk on [0,N]

3.1 *Expected number of arrivals and probability of arrival*

To obtain the expected number of arrivals, we define:

$$v_i = p_0 - (1 - r_0)\xi_i \quad (i = 1,2)$$
$$w_i = q_N - (1 - r_N)\xi_i^{-1} \quad (i = 1,2)$$

**Theorem 3** The expected number of arrivals in a [pqrs] random walk on [0,N] when starting in $i_0$ is:

$$x_0 = \frac{(\xi_1^{N-1-i_0} w_2 - \xi_2^{N-1-i_0} w_1)}{(\xi_2^{N-2} v_2 w_1 - \xi_1^{N-2} v_1 w_2)}$$

$$x_n = \frac{\zeta_1(\xi_2^{N-1-i_0} w_1 - \xi_1^{N-1-i_0} w_2)(\xi_1^{n-1} v_1 - \xi_2^{n-1} v_2)}{(\xi_2^{N-2} v_2 w_1 - \xi_1^{N-2} v_1 w_2)} \quad (1 \le n \le i_0)$$

$$x_n = \frac{\zeta_1(\xi_1^{n-1}\xi_2^{N-2} w_1 - \xi_1^{N-2}\xi_2^{n-1} w_2)(\xi_2^{1-i_0} v_1 - \xi_1^{1-i_0} v_2)}{(\xi_2^{N-2} v_2 w_1 - \xi_1^{N-2} v_1 w_2)} \quad (i_0 \le n \le N-1)$$

$$x_N = \frac{(p/q)^{N-2}(\xi_2^{1-i_0} v_1 - \xi_1^{1-i_0} v_2)}{(\xi_2^{N-2} v_2 w_1 - \xi_1^{N-2} v_1 w_2)}$$

Proof

We start in $i_0$ ($2 \le i_0 \le N-2$). We have the set of difference equations:
$$x_n = \delta(n, i_0) + p x_{n-1} + q x_{n+1} + r x_n \quad (2 \le n \le N-2)$$
with general solution, using Theorem 1:
$$x_n = \begin{cases} \zeta_1 \xi_1^{n-i_0} + C_1 \xi_1^n + C_2 \xi_2^n & (1 \le n \le i_0) \\ \zeta_1 \xi_2^{n-i_0} + C_1 \xi_1^n + C_2 \xi_2^n & (i_0 \le n \le N-1) \end{cases}$$

Notice that the solution is in this case also valid for n=1 and n=N-1, which can be seen by observing the n=2 and n=N-1 equations: $x_1$ and $x_{N-1}$ satisfy the difference pattern.

We now have:
$$x_0 = q x_1 + r_0 x_0 \Rightarrow (1 - r_0) x_0 = q(\zeta_1 \xi_1^{1-i_0} + C_1 \xi_1 + C_2 \xi_2)$$
$$x_1 = p_0 x_0 + q x_2 + r x_1 \Rightarrow p_0 x_0 = p(\zeta_1 \xi_1^{-i_0} + C_1 + C_2)$$
$$x_{N-1} = p x_{N-2} + q_N x_N + r x_{N-1} \Rightarrow q_N x_N = q(\zeta_1 \xi_2^{N-i_0} + C_1 \xi_1^N + C_2 \xi_2^N)$$
$$x_N = p x_{N-1} + r_N x_N \Rightarrow (1 - r_N) x_N = p(\zeta_1 \xi_2^{N-1-i_0} + C_1 \xi_1^{N-1} + C_2 \xi_2^{N-1})$$

We solve the four equations and get:
$$C_1 = \frac{\zeta_1 \xi_2^{N-2} w_1 \left(\xi_2^{1-i_0} v_1 - \xi_1^{1-i_0} v_2\right)}{\xi_1 (\xi_2^{N-2} v_2 w_1 - \xi_1^{N-2} v_1 w_2)}$$

$$C_2 = \frac{\zeta_1 v_2 \left(\xi_1^{N-1-i_0} w_2 - \xi_2^{N-1-i_0} w_1\right)}{\xi_2 (\xi_2^{N-2} v_2 w_1 - \xi_1^{N-2} v_1 w_2)}$$

which gives the desired result.

After some calculation we find that our result is also valid for $i_0 = 0, 1, N-1, N$.

We now can obtain directly the probability of arrival:

$$f_{ij} = \frac{x_{ij}}{x_{jj}} = \begin{cases} \dfrac{\xi_2^{1-i} v_1 - \xi_1^{1-i} v_2}{\xi_2^{1-j} v_1 - \xi_1^{1-j} v_2} & (0 \leq i < j \leq N) \\ \dfrac{\xi_2^{N-1-i} w_1 - \xi_1^{N-1-i} w_2}{\xi_2^{N-1-j} w_1 - \xi_1^{N-1-j} w_2} & (0 \leq j < i \leq N) \end{cases}$$

$$f_{ii} = 1 - x_{ii}^{-1} = \begin{cases} 1 - \dfrac{(\xi_2^{N-2} v_2 w_1 - \xi_1^{N-2} v_1 w_2)}{(\xi_1^{N-1} w_2 - \xi_2^{N-1} w_1)} & (i = 0) \\ 1 - \dfrac{(\xi_2^{N-2} v_2 w_1 - \xi_1^{N-2} v_1 w_2)}{\zeta_1 (\xi_2^{N-1-i} w_1 - \xi_1^{N-1-i} w_2)(\xi_1^{i-1} v_1 - \xi_2^{i-1} v_2)} & (0 < i < N) \\ 1 - \dfrac{(\xi_2^{N-2} v_2 w_1 - \xi_1^{N-2} v_1 w_2)}{(p/q)^{N-2} (\xi_2^{1-N} v_1 - \xi_1^{1-N} v_2)} & (i = N) \end{cases}$$

Absorption probabilities are now given by:
$P(\text{absorption in } i) = s_i x_i$ $(i = 0, ..., N)$, where $s_i = s$ $(i = 1, ..., N-1)$

3.2 *Expected time before absorption*
**Theorem 4** The expected time before absorption in i (i=0,1,....N) in a [pqrs] random walk on [0,N] is

$$m_i = \frac{1-s}{s} + \frac{(1 - \frac{s_0}{s})[w_1 \xi_2^{N-i-1} - w_2 \xi_1^{N-i-1}] + (1 - \frac{s_N}{s})(\frac{p}{q})^{N-2}[v_1 \xi_2^{1-i} - v_2 \xi_1^{1-i}]}{v_1 w_2 \xi_1^{N-2} - v_2 w_1 \xi_2^{N-2}}$$

Proof:
$$(1 - r_0) m_0 = p_0 m_1 + 1 - s_0$$
$$(1 - r) m_i = p m_{i+1} + q m_{i-1} + 1 - s \quad (i = 1, 2, ...... N-1)$$
$$(1 - r_N) m_N = q_N m_{N-1} + 1 - s_N$$

We use theorem 2 to get the desired result.

**4. [pqrs] Random walk on** $[0, \infty)$

We have probabilities $p_0, r_0, s_0$ $(p_0 + r_0 + s_0 = 1, p_0 s_0 > 0)$ for moving forward, staying for a moment and absorption in state 0. On all positive integers we have a [pqrs] walk.

4.1 *Expected number of arrivals and probability of arrival*
We start in $i_0$ with $2 \leq i_0 \leq N-2$.
We have the set of difference equations:
$$x_n = \delta(n, i_0) + p x_{n-1} + q x_{n+1} + r x_n \quad (n \geq 2)$$
with general solution:
$$x_n = \begin{cases} \zeta_1 \xi_1^{n-i_0} + C_2 \xi_2^n & (1 \leq n \leq i_0) \\ \zeta_1 \xi_2^{n-i_0} + C_2 \xi_2^n & (n \geq i_0) \end{cases}$$

Notice that the solution of theorem 1 is in this case also valid for n=1, which can be seen by observing the n=2 equation: $x_1$ satisfies the difference pattern.

We also have:
$$x_0 = qx_1 + r_0 x_0 \Rightarrow (1-r_0)x_0 = q(\zeta_1 \xi_1^{1-i_0} + C_2 \xi_2)$$
$$x_1 = p_0 x_0 + qx_2 + rx_1 \Rightarrow p_0 x_0 = p(\zeta_1 \xi_1^{-i_0} + C_2)$$

Solving this two equations leads to:

**Theorem 5** The expected number of arrivals in a [pqrs] random walk on $[0,\infty)$ is
$$x_0 = -v_1^{-1} \xi_1^{1-i_0}$$
$$x_n = \begin{cases} \zeta_1 \xi_1^{1-i_0}(\xi_1^{n-1} - \xi_2^{n-1} v_2 v_1^{-1}) & (0 < n \leq i_0) \\ \zeta_1 \xi_2^{n-1}(\xi_2^{1-i_0} - \xi_1^{1-i_0} v_2 v_1^{-1}) & (n \geq i_0) \end{cases}$$

We achieve the same result when taking $N \to \infty$ in theorem 3.

We also get:
$$f_{ij} = \frac{x_{ij}}{x_{jj}} = \begin{cases} \xi_1^{j-i} & (0 \leq j < i) \\ \dfrac{\xi_2^{1-i} v_1 - \xi_1^{1-i} v_2}{\xi_2^{1-j} v_1 - \xi_1^{1-j} v_2} & (0 \leq i < j) \end{cases}$$

$$f_{ii} = 1 - x_{ii}^{-1} = \begin{cases} 1 + \dfrac{v_1}{\xi_1} & (i=0) \\ 1 - \dfrac{\xi_1^{i-1} v_1}{\zeta_1(\xi_1^{i-1} v_1 - \xi_2^{i-1} v_2)} & (i>0) \end{cases}$$

After some calculation we find that the results are also valid for $i_0 = 0, 1, N-1, N$.

4.2 *Expected time before absorption*

**Theorem 6** The expected time before absorption in a [pqrs] random walk on $[0,\infty)$ when starting in i is:
$$m_i = \frac{1-s}{s} + v_1^{-1}(1 - \frac{s_0}{s})\xi_1^{1-i} \quad (i = 0,1,2,.....)$$

Proof:
$$(1-r_0)m_0 = p_0 m_1 + 1 - s_0$$
$$(1-r)m_i = pm_{i+1} + qm_{i-1} + 1 - s \quad (i > 0)$$

Use $m_i = a\xi_1^{-i} + \dfrac{1-s}{s}$ to get the desired result, or take $N \to \infty$ in theorem 4.

4.3 *Expected time before absorption in j, when starting in i.*

**Theorem 7** The expected time before absorption in j in a [pqrs] random walk on $[0,\infty)$ when starting in i is

If j=0 then
$$m_{i0} = s\zeta_1 v_1^{-2} \xi_1^{1-i}(p_0 - iv_1)$$

If $1 \leq j \leq i$ then $m_{ij} =$
$$s\zeta_1^2 v_1^{-1} \xi_1^{1-i}\{[\alpha\zeta_1 - (1-r_0)v_1^{-1}\xi_1 + i - 1](v_1\xi_1^{j-1} - v_2\xi_2^{j-1}) + (1-r_0)(\xi_1^j + \xi_2^j) - (j-1)(v_1\xi_1^{j-1} + v_2\xi_2^{j-1})\}$$

If $j \geq i$ then $m_{ij} =$
$$s\zeta_1^2 v_1^{-1} \xi_2^{j-1}\{[\alpha\zeta_1 - (1-r_0)v_1^{-1}\xi_1 + j - 1](v_1\xi_2^{1-i} - v_2\xi_1^{1-i}) + p(1-r_0)(\xi_1^{-i} + \xi_2^{-i}) - (i-1)(v_1\xi_2^{1-i} + v_2\xi_1^{1-i})\}$$

Proof:
We have:

$$m_{ij} = s\sum_{k=0}^{\infty} kp_{ij}^{(k)} = s\left(\frac{dX_{ij}}{dz}\right)_{z=1} \quad (i,j \geq 0)$$

where (extension of theorem 5):

$$X_{i0}(z) = -v_1^{-1}\xi_1^{1-i}$$

$$X_{ij}(z) = \begin{cases} \zeta_z \xi_1^{1-i}(\xi_1^{n-1} - \xi_2^{n-1}v_2v_1^{-1}) & (0 < j \leq i) \\ \zeta_z \xi_2^{n-1}(\xi_2^{1-i} - \xi_1^{1-i}v_2v_1^{-1}) & (j \geq i) \end{cases}$$

Implicit differentiation of $qz\xi^2 - (1-rz)\xi + pz = 0$ gives:

$$\left(\frac{d\xi_i}{dz}\right)_{z=1} = (-1)^i \zeta_1 \xi_i \quad (i=1,2)$$

We also have:

$$\zeta_z = [(1-rz)^2 - 4pqz^2]^{-\frac{1}{2}} \Rightarrow \left(\frac{d\zeta}{dz}\right)_{z=1} = \zeta_1^3[r(1-r) + 4pq]$$

and after some calculation we find the desired result.

## 5. [pqrs] Random walk on $(-\infty, \infty)$

5.1 *Expected number of arrivals and probability of arrival*

**Theorem 8** The expected number of arrivals in state n, when starting in $i_0$ in a [pqrs] random walk on $(-\infty, \infty)$ is

$$x_n = \begin{cases} \zeta_1 \xi_1^{n-i_0} & (n \leq i_0) \\ \zeta_1 \xi_2^{n-i_0} & (n \geq i_0) \end{cases}$$

The probability of arrival in state j when starting in i a [pqrs] random walk on $(-\infty, \infty)$ is

$$f_{ij} = \frac{x_{ij}}{x_{jj}} = \begin{cases} \xi_1^{j-i} & (j < i) \\ \xi_2^{j-i} & (j > i) \end{cases}$$

$$f_{ii} = \frac{x_{ii}-1}{x_{ii}} = 1 - \sqrt{(1-r)^2 - 4pq}$$

Proof
Starting in $i_0$, we have the set of difference equations:

$$x_n = \delta(n, i_0) + px_{n-1} + qx_{n+1} + rx_n \quad (n \in \mathbb{Z})$$

Using Theorem 1 and $C_1 = C_2 = 0$ (because of $\xi_1 > 1$ and $0 < \xi_2 < 1$) gives the desired result.

5.2 *Expected time before absorption and expected time before absorption in j.*

**Theorem 9** The expected time before absorption in a [pqrs] random walk on $(-\infty, \infty)$ is

$$m = \frac{1-s}{s}$$

The expected time before absorption in j when starting in 0 in a [pqrs] random walk on $(-\infty, \infty)$ is

$$m_{0j} = \begin{cases} s\zeta_1^2 \xi_1^j \{\zeta_1[r(1-r)+4pq] - j\} & (j \leq 0) \\ s\zeta_1^2 \xi_2^j \{\zeta_1[r(1-r)+4pq] + j\} & (j \geq 0) \end{cases}$$

Proof

The expected time before absorption $m$ satisfies: $m = p(m+1) + q(m+1) + r(m+1)$, so: $m = \dfrac{1-s}{s}$

We now study the expected time for absorption in j, when starting in 0 (this is no limitation in the $(-\infty, \infty)$ case).

$$m_{0j} = s\sum_{k=0}^{\infty} k p_{0j}^{(k)} = s\left(\dfrac{dX_{0j}}{dz}\right)_{z=1}$$

where
$$X_{0j} = \begin{cases} \zeta_z \xi_1^j(z) & (j \leq 0) \\ \zeta_z \xi_2^j(z) & (j \geq 0) \end{cases}$$

We find the desired result by following the proof of theorem 7.

5.3 *Formulae for n step probability.*

**Theorem 10** The n step probability for a [pqrs] random walk on $(-\infty, \infty)$ when starting in 0 is

$$p_{0k}^{(n)} = \begin{cases} p^{\frac{n+k}{2}} q^{\frac{n-k}{2}} \sum_{w=0}^{\frac{n-k}{2}} \left[\dfrac{r^2}{pq}\right]^w \binom{n}{2w}\binom{n-2w}{\frac{n-k}{2}-w} & \text{if } n = k \pmod 2 \\ \\ p^{\frac{n+k-1}{2}} q^{\frac{n-k-1}{2}} r \sum_{w=0}^{\frac{n-k-1}{2}} \left[\dfrac{r^2}{pq}\right]^w \binom{n}{2w+1}\binom{n-2w-1}{\frac{n-k-1}{2}-w} & \text{if } n = k+1 \pmod 2 \end{cases}$$

we also have

$$p_{0k}^{(n)} = q^n (r + \sqrt{r^2 - 4pq})^{-k-n} \sum_{m=k}^{n} \binom{n}{m}\binom{n}{m-k}\left(\dfrac{-r - \sqrt{r^2 - 4pq}}{-r + \sqrt{r^2 - 4pq}}\right)^m \quad (-n \leq k \leq n)$$

Proof
To obtain formulae for n step probability we use two methods.
*Combinatorial approach*
We have

$$p_{0k}^{(n)} = \sum \binom{n}{n_p, n_q, n_r} p^{n_p} q^{n_q} r^{n_r} \quad \text{with} \quad n_p + n_q + n_r = n, \quad n_p - n_q = k$$

Case n=k (mod 2)
By taking $n_p = \dfrac{n+k}{2} - w$, $n_q = \dfrac{n-k}{2} - w$, $n_r = 2w$ $(w = 0,1..., \dfrac{n-k}{2})$ we obtain the first result.

Case n=k+1 (mod 2)
By taking $n_p = \dfrac{n+k-1}{2} - w$, $n_q = \dfrac{n-k-1}{2} - w$, $n_r = 2w+1$ $(w = 0,1..., \dfrac{n-k-1}{2})$ we obtain the next result.

*PGF approach*

We define a generating function $\quad f(z) = pz + \dfrac{q}{z} + r$

So:  $[f(z)]^n = \left[pz + \dfrac{q}{z} + r\right]^n = z^{-n}[pz^2 + rz + q]^n = p^n z^{-n}(z - z_1)^n (z - z_2)^n =$

$$q^n z^{-n}(1 - \dfrac{z}{z_1})^n (1 - \dfrac{z}{z_2})^n = q^n z^{-n} \left[\sum_{m=0}^{n} \binom{n}{m}(\dfrac{-z}{z_1})^m\right]\left[\sum_{l=0}^{n} \binom{n}{l}(\dfrac{-z}{z_2})^l\right].$$

Focussing on the term with $z^k$ (using $m + l = k + n$):

$$p_{0k}^{(n)} = q^n z^{-n} \left[\sum_{m=k}^{n} \binom{n}{m}(\dfrac{-z}{z_1})^m \binom{n}{n+k-m}(\dfrac{-z}{z_2})^{k+n-m}\right] =$$

$$q^n (-z_2)^{-k-n} \sum_{m=k}^{n} \binom{n}{m}\binom{n}{m-k}(\dfrac{z_2}{z_1})^m \quad (-n \leq k \leq n)$$

where
$$z_{1,2} = \dfrac{-r \pm \sqrt{r^2 - 4pq}}{2p}$$

## 6. Modified [pqrs] random walk on [0,N].

We consider a random walk on the interval [0,N] with the following specifications:
between 0 and M (M<N) there is a $[p_2, q_2, r_2, s_2]$ walk, between M and N we have a
$[p_1, q_1, r_1, s_1]$ walk and in 0,M and N we have specific probabilities $[p_i, q_i, r_i, s_i]$ (i=0,M,N), where
$p_N = q_0 = 0$. We use the shortcut M[pqrs] for this modified [pqrs] random walk.
We start in $i_0$. We handle the case $0 < M < i_0 < N$ ($0 < i_0 < M < N$ can be handled the same way or
direct by the result of this section: use a reflection in 0, followed by a translation).

6.1 *Expected number of arrivals and probability of arrival*
**Theorem 11** The expected number of arrivals in n in a M[pqrs] random walk on [0,N] when starting
in $i_0$ is:

$$\dfrac{p_2 q_M \phi_0}{p_0 q_2 \phi_M} x_M \quad (n = 0)$$

$$\dfrac{q_M \phi_n}{q_2 \phi_M} x_M \quad (1 \leq n \leq M-1)$$

$$\dfrac{\lambda_1 \left[\dfrac{q_1}{p_1}\right]^{i_0}}{(\gamma_1 \mu_2 - \gamma_2 \mu_1)} \quad (n = M)$$

$$\dfrac{\left[\dfrac{q_1}{p_1}\right]^{i_0}(\mu_1 \xi_1^n - \mu_2 \xi_2^n)}{(\gamma_1 \mu_2 - \gamma_2 \mu_1)} \quad (M+1 \leq n \leq i_0)$$

$$\beta_n x_{i_0} \quad (i_0 \leq n \leq N-1)$$

$$\dfrac{q_1}{q_N} \beta_N x_{i_0} \quad (n = N)$$

where
$$\lambda_i = \sqrt{(1 - r_i)^2 - 4 p_i q_i} \quad (i = 1, 2)$$

$$\eta_i = \frac{(1-r_2) \pm \lambda_2}{2q_2} \; ; \quad \xi_i = \frac{(1-r_1) \pm \lambda_1}{2q_1} \quad (i=1,2)$$

$$t_i = p_2(1-r_0) - p_0 q_2 \eta_i \quad (i=1,2); \quad t = \frac{t_1}{t_2}$$

$$u_i = \xi_i^{N-1}[q_1(1-r_N)\xi_i - p_1 q_N] \quad (i=1,2); \quad u = \frac{u_1}{u_2}$$

$$\varphi_n = \eta_1^n - t\eta_2^n; \quad \beta_n = \frac{\xi_1^n - u\xi_2^n}{\xi_1^{i_0} - u\xi_2^{i_0}}$$

$$\gamma_i = \xi_i^{-i_0}[r_1 - 1 + q_1(\beta_{i_0+1} + \xi_i)] \quad (i=1,2)$$

$$\mu_i = \xi_i^{-M}[1 - r_M - p_M \xi_i^{-1} - \frac{p_2 q_M \varphi_{M-1}}{q_2 \varphi_M})] \quad (i=1,2)$$

Proof
We have:
$$x_n = p_2 x_{n-1} + q_2 x_{n+1} + r_2 x_n \quad (2 \le n \le M-2)$$
$$x_n = p_1 x_{n-1} + q_1 x_1 + r_1 x_n + \delta(n,i_0) \quad (M+2 \le n \le N-2)$$

We note that both $x_1, x_{M-1}, x_{M+1}$ and $x_{N-1}$ are present in one of these formulae and therefore satisfies the relevant pattern:

$$x_n = a\eta_1^n + b\eta_2^n \quad (1 \le n \le M-1), \text{ where } q_2\eta^2 - (1-r_2)\eta + p_2 = 0$$
$$x_n = c\xi_1^n + d\xi_2^n \quad (M+1 \le n \le i_0), \text{ where } q_1\xi^2 - (1-r_1)\xi + p_1 = 0$$
$$x_n = h_1\xi_1^n + h_2\xi_2^n \quad (i_0 \le n \le N-1)$$

We now have:
$$x_0 = q_2 x_1 + r_0 x_0 \Rightarrow (1-r_0)x_0 = q_2(a\eta_1 + b\eta_2)$$
$$x_1 = p_0 x_0 + q_2 x_2 + r_2 x_1 \Rightarrow p_0 x_0 = p_2(a+b)$$
$$x_{M-1} = p_2 x_{M-2} + q_M x_M + r_2 x_{M-1} \Rightarrow q_M x_M = q_2(a\eta_1^M + b\eta_2^M)$$
$$x_M = p_2 x_{M-1} + q_1 x_{M+1} + r_M x_M \Rightarrow (1-r_M)x_M = p_2(a\eta_1^{M-1} + b\eta_2^{M-1}) + q_1(c\xi_1^{M+1} + d\xi_2^{M+1})$$
$$x_{M+1} = p_M x_M + q_1 x_{M+2} + r_1 x_{M+1} \Rightarrow p_M x_M = p_1(c\xi_1^M + d\xi_2^M)$$
$$x_{i_0} = c\xi_1^{i_0} + d\xi_2^{i_0}$$
$$(1-r_1)x_{i_0} = 1 + p_1(c\xi_1^{i_0-1} + d\xi_2^{i_0-1}) + q_1(h_1\xi_1^{i_0+1} + h_2\xi_2^{i_0+1})$$
$$x_{i_0} = h_1\xi_1^{i_0} + h_2\xi_2^{i_0}$$
$$x_{N-1} = p_1 x_{N-2} + q_N x_N + r_1 x_{N-1} \Rightarrow q_N x_N = q_1(h_1\xi_1^N + h_2\xi_2^N)$$
$$x_N = p_1 x_{N-1} + r_N x_N \Rightarrow (1-r_N)x_N = p_1(h_1\xi_1^{N-1} + h_2\xi_2^{N-1})$$

Solving this 10 equations gives the desired result.

*6.1.1 Special case 1: asymmetric random walk with three mfb's.*

When we choose $p_1 = p_2 = p$; $q_1 = q_2 = q$; $r_1 = r_2 = s_1 = s_2 = 0$; $p \ne q$, we have an asymmetric random walk with three mfb's in 0,M and N.

We can use theorem 11, where $\eta_1 = \xi_1 = \frac{p}{q}$; $\quad \eta_2 = \xi_2 = 1$; $\quad \lambda_1 = \lambda_2 = |p-q|$

*6.1.2 Special case 2: symmetric random walk with three mfb's.*

Taking $p_1 = p_2 = q_1 = q_2 = \frac{1}{2}$ we get a simple symmetric random walk with mfb's in 0,M,N.

Expected number of arrivals are now given by:

$$x_n = an + b \quad (1 \le n \le M - 1)$$
$$x_n = cn + d \quad (M + 1 \le n \le i_0)$$
$$x_n = h_1 n + h_2 \quad (i_0 \le n \le N - 1)$$

Proceeding along the same lines as in theorem 11 we get the expected number of arrivals in n:

$$\frac{q_M (\Omega - i_0)}{(p_0 + M s_0) k(\Omega)} \quad (n = 0)$$

$$\frac{2 q_M (n + \frac{p_0}{s_0})(\Omega - i_0)}{(M + \frac{p_0}{s_0}) k(\Omega)} \quad (1 \le n \le M - 1)$$

$$\frac{(\Omega - i_0)}{k(\Omega)} \quad (n = M)$$

$$\frac{2(\Omega - i_0) k(n)}{k(\Omega)} \quad (M + 1 \le n \le i_0)$$

$$\frac{2(\Omega - n) k(i_0)}{k(\Omega)} \quad (i_0 \le n \le N - 1)$$

$$\frac{k(i_0)}{s_N k(\Omega)} \quad (n = N)$$

where

$$\Omega = N + \frac{q_N}{s_N}; \quad k(z) = p_M + \theta(z - M) \; ; \; \theta = s_M + \frac{q_M}{(M + \frac{p_0}{s_0})}$$

We get the same result by applying de l'Hospitals rule in theorem 11.
It is easily verified that P(absorption)= $s_0 x_0 + s_M x_M + s_N x_N = 1$

## 7. Modified [pqrs] random walk on $[0, \infty)$

*7.1 Expected number of arrivals and probability of arrival*
We consider a random walk on the interval $[0, \infty)$ with the following specifications:
between 0 and M there is a $[p_2, q_2, r_2, s_2]$ walk, beyond M we have a $[p_1, q_1, r_1, s_1]$ walk and in 0 and M we have specific probabilities $[p_i, q_i, r_i, s_i]$ (i=0,M), where $q_0 = 0$.
We start in $i_0$ and we handle the case $0 < M < i_0$.

We suppose $p_i q_i s_i > 0$ $(i = 1,2)$, which implies $\xi_1 = \frac{(1 - r_1) + \lambda_1}{2 q_1} > 1$ and $0 < \xi_2 = \frac{(1 - r_1) - \lambda_1}{2 q_1} < 1$

We now have: $x_n = h_2 \xi_2^n \quad (n \ge i_0)$
We can now use the same technique as in theorem 11, now solving 9 equations.
We get the same result by taking the solution on [0,N] and letting $N \to \infty$, which results in:
$u \to \infty$ and $\beta_n \to \xi_2^{n - i_0}$.

*7.2 Expected time before absorption*

**Theorem 12** The expected time before absorption in a M[pqrs] random walk on $[0, \infty)$ when starting in i is
If $0 \le i \le M$ then:

$$m_i = \sigma_2 + \{q[1-(1+\sigma_2)s_0](\frac{q}{p})^{i-1}[\alpha_M(\xi_1^{i-M} - \xi_2^{i-M}) - q_M(\xi_1^{i-M+1} - \xi_2^{i-M+1})] +$$

$$q[1-(1+\sigma_2)s_M + p_M(\sigma_1 - \sigma_2)(1-\xi_1^M)](v_1\xi_2^{1-i} - v_2\xi_1^{1-i})\} *$$

$$\{p(1-r_0)\alpha_M(\xi_1^{-M} - \xi_2^{-M}) + q[p_0\alpha_M + (1-r_0)q_M](\xi_1^{1-M} - \xi_2^{1-M}) + qp_0q_M(\xi_1^{2-M} - \xi_2^{2-M})\}^{-1}$$

If $i \geq M$ then:

$$m_i = \sigma_1 + (m_M - \sigma_1)\xi_1^{M-i}$$

where

$$\sigma_k = \frac{1-s_k}{s_k} \quad (k=1,2)$$

$$\alpha_M = 1 - r_M - p_M \xi_1^M$$

$$v_i = p_0 - (1-r_0)\xi_i \quad (i=1,2)$$

Proof
The expected time before absorption when starting in i satisfies:

$$(1-r_2)m_i = p_2 m_{i+1} + q_2 m_{i-1} + 1 - s_2 \quad (1 \leq i \leq M-1)$$
$$(1-r_1)m_i = p_1 m_{i+1} + q_1 m_{i-1} + 1 - s_1 \quad (i \geq M+1)$$

with solution:

$$m_i = \begin{cases} a_2\xi_1^{-i} + b_2\xi_2^{-i} + \frac{1-s_2}{s_2} & (0 \leq i \leq M) \\ a_1\xi_1^{-i} + \frac{1-s_1}{s_1} & (i \geq M) \end{cases}$$

We also have:

$$(1-r_0)m_0 = p_0 m_1 + 1 - s_0$$
$$m_0 = a_2 + b_2 + \sigma_2$$
$$(1-r_M)m_M = p_M m_{M+1} + 1 - s_M$$
$$m_M = a_1 \xi_1^{-M} + \sigma_1$$
$$m_M = a_2 \xi_1^{-M} + b_2 \xi_2^{-M} + \sigma_2$$

After some calculations we get our result.

## 8. Modified [pqrs] random walk on $(-\infty, \infty)$

8.1 *Expected number of arrivals and probability of arrival*
We consider a random walk on the interval $(-\infty, \infty)$ with the following specifications:
On the negative integers there is a $[p_2, q_2, r_2, s_2]$ walk, on the positive integers we have a $[p_1, q_1, r_1, s_1]$ walk and in 0 we have specific probabilities $[p_0, q_0, r_0, s_0]$.
We start in $i_0$ and we handle the case $i_0 > 0$.
We suppose $p_i q_i s_i > 0$ $(i=1,2)$, which implies $\xi_1 > 1$, $0 < \xi_2 < 1$, $\eta_1 > 1$, $0 < \eta_2 < 1$.
We have:

$$x_n = p_2 x_{n-1} + q_2 x_{n+1} + r_2 x_n \quad (n \leq -2)$$
$$x_n = p_1 x_{n-1} + q_1 x_{n+1} + r_1 x_n + \delta(n, i_0) \quad (n \geq 2)$$

Solutions are given by:

$$x_n = a_1\eta_1^n \quad (n \leq -1)$$
$$x_n = c\xi_1^n + d\xi_2^n \quad (1 \leq n \leq i_0)$$
$$x_n = h_2\xi_2^n \quad (n \geq i_0)$$

We now have:

$$x_{-1} = p_2 x_{-2} + q_0 x_0 + r_2 x_{-1} \Rightarrow q_0 x_0 = q_2 a_1$$
$$x_0 = p_2 x_{-1} + q_1 x_1 + r_0 x_0 \Rightarrow (1-r_0)x_0 = p_2 a_1 \eta_1^{-1} + q_1(c\xi_1 + d\xi_2)$$
$$x_1 = p_0 x_0 + q_1 x_2 + r_1 x_1 \Rightarrow p_0 x_0 = p_1(c+d)$$
$$x_{i_0} = c\xi_1^{i_0} + d\xi_2^{i_0}$$
$$(1-r_1)x_{i_0} = 1 + p_1(c\xi_1^{i_0-1} + d\xi_2^{i_0-1}) + q_1 h_2 \xi_2^{i_0+1}$$
$$x_{i_0} = h_2 \xi_2^{i_0}$$

Solving this six equations gives:

**Theorem 13** The expected number of arrivals in state n, when starting in $i_0$ in a M[pqrs] random walk on $(-\infty,\infty)$ is

$$x_n = \begin{cases} \dfrac{q_0 \eta_1^n}{q_2 \tau_{1,i_0}} & (n<0) \\ \tau_{1,i_0}^{-1} & (n=0) \\ \lambda_1^{-1}[\xi_1^{n-i_0} - \dfrac{\tau_{2,n}}{\tau_{1,i_0}}] & (0<n \leq i_0) \\ \lambda_1^{-1}\xi_2^{n-i_0}[1 - \dfrac{\tau_{2,i_0}}{\tau_{1,i_0}}] & (n \geq i_0) \end{cases}$$

where

$$\tau_{i,n} = (1 - r_0 - p_2 q_0 q_2^{-1}\eta_1^{-1} - p_0\xi_i^{-1})\xi_i^n \quad (i=1,2)$$

Remark. If $i_0 = 0$ we get the simpler result:

$$x_n = \begin{cases} \dfrac{q_0 \eta_1^n}{q_2 \tau_{1,0}} & (n<0) \\ \tau_{1,0}^{-1} & (n=0) \\ \dfrac{p_0 \xi_2^n}{p_1 \tau_{1,0}} & (n>0) \end{cases}$$

*8.1.1. Probability of absorption in a special case.*
When we choose $p_1 + q_1 + r_1 = 1$; $p_2 + q_2 + r_2 = 1$; $s_1 = s_2 = 0$, we have an asymmetric classic random walk with different parameters left and right of the origin.
We first consider the case $p_1 > q_1$ and $p_2 > q_2$ and we find directly that the probability of absorption is:

$$s_0 x_0 = \dfrac{s_0}{[s_0 + p_0(1 - \dfrac{q_1}{p_1})]\left[\dfrac{p_1}{q_1}\right]^{i_0}}$$

$$P(\text{escape in } \infty) = (p_1 - q_1)x_{i_0} = 1 - \frac{\tau_{2,i_0}}{\tau_{1,i_0}} = 1 - \frac{s_0}{[s_0 + p_0(1 - \frac{q_1}{p_1})]\left[\frac{p_1}{q_1}\right]^{i_0}}$$

If $p_1 > q_1$ and $p_2 < q_2$ we get:

$$P(\text{absorption in } 0) = \frac{s_0}{[s_0 + p_0(1 - \frac{q_1}{p_1}) + q_0(1 - \frac{p_2}{q_2})]\left[\frac{p_1}{q_1}\right]^{i_0}}$$

$$P(\text{escape in } \infty) = 1 - \frac{s_0 + q_0(1 - \frac{p_2}{q_2})}{[s_0 + p_0(1 - \frac{q_1}{p_1}) + q_0(1 - \frac{p_2}{q_2})]\left[\frac{p_1}{q_1}\right]^{i_0}}$$

$$P(\text{escape in } -\infty) = \frac{q_0(1 - \frac{p_2}{q_2})}{[s_0 + p_0(1 - \frac{q_1}{p_1}) + q_0(1 - \frac{p_2}{q_2})]\left[\frac{p_1}{q_1}\right]^{i_0}}$$

8.2 *Expected time before absorption*

We have:

$$(1 - r_2)m_i = p_2 m_{i+1} + q_2 m_{i-1} + 1 - s_2 \quad (i \leq -1)$$
$$(1 - r_1)m_i = p_1 m_{i+1} + q_1 m_{i-1} + 1 - s_1 \quad (i \geq 1)$$

Solution:

$$m_i = \begin{cases} b\xi_2^{-i} + \frac{1 - s_2}{s_2} & (i \leq 0) \\ a\xi_1^{-i} + \frac{1 - s_1}{s_1} & (i \geq 0) \end{cases}$$

So:

$$m_0 = b + \frac{1 - s_2}{s_2} = a + \frac{1 - s_1}{s_1}$$

We also have:

$$(1 - r_0)m_0 = p_0 m_1 + q_0 m_{-1} + 1 - s_0$$

Solving this three equations leads to:

**Theorem 14** The expected time before absorption in a M[pqrs] random walk on $(-\infty, \infty)$ when starting in i is

$$m_i = \begin{cases} \dfrac{[(\sigma_2 - \sigma_1)p_0(1-\xi_1^{-1}) - \sigma_2 s_0]}{[s_0 + p_0(1-\xi_1^{-1}) + q_0(1-\xi_2)]} \xi_2^{-i} + \sigma_2 & (i \leq 0) \\ \dfrac{[(\sigma_2 - \sigma_1)q_0(1-\xi_2) - \sigma_1 s_0]}{[s_0 + p_0(1-\xi_1^{-1}) + q_0(1-\xi_2)]} \xi_1^{-i} + \sigma_1 & (i \geq 0) \end{cases}$$

where $\sigma_k = \dfrac{1-s_k}{s_k}$ $(k=1,2)$